\documentclass{article}
\begin{document}
\newtheorem{proposition}{Proposition}[section]
\newtheorem{definition}{Definition}[section]
\newtheorem{lemma}{Lemma}[section]
\newcommand{\xl}{\stackrel{\rightharpoonup}{\cdot}}
\newcommand{\xr}{\stackrel{\leftharpoonup}{\cdot}}
\newcommand{\xlplus}{\stackrel{\rightharpoonup}{+}}
\newcommand{\xrplus}{\stackrel{\leftharpoonup}{+}}
\newcommand{\xluplus}{\stackrel{\rightharpoonup}{\uplus}}
\newcommand{\xruplus}{\stackrel{\leftharpoonup}{\uplus}}
\newcommand{\xlodot}{\stackrel{\rightharpoonup}{\odot}}
\newcommand{\xrodot}{\stackrel{\leftharpoonup}{\odot}}
\newcommand{\Ll}{\stackrel{\rightharpoonup}{L}}
\newcommand{\Lr}{\stackrel{\leftharpoonup}{L}}
\newcommand{\Rl}{\stackrel{\rightharpoonup}{R}}
\newcommand{\Rr}{\stackrel{\leftharpoonup}{R}}

\title{\bf A Class of Ring-like Objects}
\author{Keqin Liu\\Department of Mathematics\\The University of British Columbia
\\Vancouver, BC\\
Canada, V6T 1Z2}
\date{November 23, 2004}
\maketitle

\begin{abstract} We introduce the notions of one-sided dirings, $3$-irreducible 
left modules, $3$-primitive left dirings, $3$-semi-primitive left dirings, 
$3$-primitive ideals and $3$-radicals. The main results consists of two parts. 
The first part establishes two external characterizations of a $3$-semi-primitive 
left diring. The second part characterizes the $3$-radical of a left diring by 
using $3$-primitive ideals.
\end{abstract}

By forgetting some structures of a $7$-tuple introduced in Chapter 4 of 
\cite{Liu2}, we get three roads of generalizing the notion of a ring $R$. The first 
one is to keep the additive group structure of $R$ and to replace the 
multiplicative monoid structure of $R$ by a dimonoid with a one-sided bar-unit. 
The second one is to replace the additive group structure of $R$ by a commutative 
digroup and to keep the multiplicative monoid structure of $R$. The third one is to 
replace the additive group structure of $R$ by a commutative digroup and to replace 
the multiplicative monoid structure of $R$ by a dimonoid with a one-sided bar-unit. 
Although we do not know how far we can go along the third road now, the first 
two roads are good enough to develop the counterpart of the basic ring theory. The
 purpose of this paper is to study the counterpart of the Jacobson radical 
for rings along the first road. 

\medskip
This paper consists of five sections. In Section 1 we introduce the notion of a 
one-sided diring and discusses its basic properties. In Section 2 we consider 
some fundamental concepts and results about a left module over a left diring. 
In Section 3 we introduce the notion of a $3$-irreducible left module and prove 
that Schur Lemma is still true for $3$-irreducible left modules over a left diring. 
In Section 4 we introduce the notions of $3$-primitive left dirings and 
$3$-semi-primitive left dirings, and establish two external characterizations of a 
$3$-semi-primitive left diring. In Section 5 we introduce the notion of the 
$3$-radical of a left diring by using the intersection of the annihilators of all 
$3$-irreducible left $R$-modules, and prove that  the $3$-radical of a left diring $R$ 
is equal to the intersection of the $3$-primitive ideals of $R$.

\section{The Notion of One-sided Dirings}

We begin this section with the definition of a one-sided diring.

\begin{definition}\label{def1.1} A nonempty set $R$ is called a {\bf left diring} 
(or {\bf right diring}) if there are three binary operations $+$, $\xl$ and 
$\xr$ on $R$ such that the following three properties hold
\begin{description}
\item[(i)] $(\, R , \, + \,)$ is an Abelian group with the identity $0$. 
\item[(ii)] $(\, R, \, \xl , \, \xr \,)$ is a dimonoid with a left bar-unit 
$e_{\ell}$ (or a right bar-unit $e_r$).\footnote{The notion of a one-sided bar-unit 
of a dimonoid was introduced in Definition 6.3 of \cite{Liu2}}
\item[(iii)] The distributive laws
$$ x\ast (y+z)=x\ast y +x\ast z, \quad (y +z)\ast x= y\ast x + z\ast x $$
hold for all $x$, $y$, $z\in R$ and $\ast\in \{\, \xl ,\, \xr \, \}$.
\end{description}
\end{definition}

A left diring or a right diring is called a {\bf one-sided diring} and is denoted 
by $(\, R, \, + , \, \xl , \, \xr \,)$.

\begin{definition}\label{def1.2} A one-sided diring 
$(\, R, \, + , \, \xl , \, \xr \,)$ is called a {\bf diring} if 
$(\, R, \, \xl , \, \xr \,)$ is a dimonoid with a bar-unit. 
\end{definition}

By forgetting the vector space structure of a dialgebra with a 
bar-unit(\cite{Loday2}), we can regard a dialgebra with a bar-unit as a diring.

\medskip
If $(\, R, \, + , \, \xl , \, \xr \,)$ is a one-sided diring, then the binary 
operations $+$, $\xl$ and $\xr$ are called the {\bf addition}, 
the {\bf left product} and the {\bf right product}, respectively. The Abelian 
group $(\, R, \, + \,)$ is called the {\bf additive group} of $R$, and the 
identity $0$ of the additive group is called the {\bf zero element} of $R$. 
If $x$ is an element of $R$, then the group inverse of $x$ in the additive 
group is denoted by $-x$. A one-sided bar-unit of the dimonoid 
$(\, R, \, \xl , \, \xr \,)$ is called a 
one-sided {\bf multiplicative bar-unit} of $R$. The left halo of the dimonoid 
$(\, R, \, \xl , \, \xr \,)$ is called the {\bf left multiplicative halo} of 
$R$ and is denoted by $\hbar _{\ell}^\times (R)$. 
The right halo of the dimonoid $(\, R, \, \xl , \, \xr \,)$ is called the 
{\bf right multiplicative halo} of $R$ and is denoted by $\hbar _{r}^\times (R)$. 
The halo of the dimonoid $(\, R, \, \xl , \, \xr \,)$ is called the 
{\bf multiplicative halo} of $R$ and is denoted by $\hbar ^\times (R)$.
Thus, we have
$$
\hbar _{\ell}^\times (R)=
\{\, \alpha \in R \, | \, \mbox{$\alpha\xr x=x$ for all $x\in R$} \,\},
$$
$$
\hbar _{r}^\times (R)=
\{\, \alpha \in R \, | \, \mbox{$x=x\xl \alpha$ for all $x\in R$} \,\}
$$
and 
$$
\hbar ^\times (R)=
\{\, \alpha \in R \, | \, \mbox{$\alpha\xr x=x=x\xl \alpha$ for all $x\in R$} \,\}.
$$

\medskip
The following example gives a left diring which is not a diring.

\bigskip
\noindent
{\bf Example} Let $H: =\{\, 0, \, a ,\, b ,\, c \,\}$ be a set of four distinct 
elements. We define three binary operations $+$, $\xl$ and $\xr$ on $H$ as follows:
$$
\begin{tabular}{|c||c|c|c|c|}\hline
$+$ &$0$&$a$&$b$&$c$\\
\hline\hline
$0$&$0$&$a$&$b$&$c$\\
\hline
$a$&$a$&$0$&$c$&$b$\\
\hline
$b$&$b$&$c$&$0$&$a$\\
\hline
$c$&$c$&$b$&$a$&$0$\\
\hline
\end{tabular}
$$
$$
\begin{tabular}{|c||c|c|c|c|}\hline
$\xl$ &$0$&$a$&$b$&$c$\\
\hline\hline
$0$&$0$&$0$&$0$&$0$\\
\hline
$a$&$0$&$b$&$b$&$0$\\
\hline
$b$&$0$&$b$&$b$&$0$\\
\hline
$c$&$0$&$0$&$0$&$0$\\
\hline
\end{tabular}
\qquad\qquad\qquad
\begin{tabular}{|c||c|c|c|c|}\hline
$\xr$ &$0$&$a$&$b$&$c$\\
\hline\hline
$0$&$0$&$0$&$0$&$0$\\
\hline
$a$&$0$&$a$&$b$&$c$\\
\hline
$b$&$0$&$a$&$b$&$c$\\
\hline
$c$&$0$&$0$&$0$&$0$\\
\hline
\end{tabular}
$$

One can check that $(\, H ,\, \xl ,\, \xr \,)$ is a left diring with 
$h_{\ell}^{\times}(H)=\{a, \, b\}$. Since $h_{r}^{\times}(H)=\emptyset$, $H$ does not 
have any bar-unit. Hence, $H$ is not a diring. 
\hfill\raisebox{1mm}{\framebox[2mm]{}}

\medskip
It is clear that if $(\, R ,\, +, \, \xl , \, \xr \,)$ is a left diring with a 
left bar-unit $e_{\ell}$ , then 
$(\, \mathring{R} , \, +, \, \stackrel{\rightharpoonup}{\circ} ,\, 
\stackrel{\leftharpoonup}{\circ} \,)$ is a right diring with a right bar-unit 
$\mathring{e}$, where $\mathring{e}:= e_{\ell}$ and the binary operations 
$\stackrel{\rightharpoonup}{\circ}$ and  $\stackrel{\leftharpoonup}{\circ}$ are 
defined by
\begin{eqnarray*}
x\stackrel{\rightharpoonup}{\circ} y :&=& y\xr x,\\
x\stackrel{\leftharpoonup}{\circ} y:&=& y\xl x,
\end{eqnarray*}
where $x$, $y\in R$. $\mathring{R}$ is called the 
{\bf opposite one-sided diring} of an one-sided diring $R$. Using the opposite 
one-sided diring, a fact about a left diring can be converted to a fact about a 
right diring, and vice versa. Hence, we will only discuss left dirings.

\medskip
The following are direct consequences of the distributive laws:
$$ x\ast 0=0=0\ast x,$$
$$ (-x)\ast y=x\ast (-y)=-(x\ast y),$$
where $x$, $y$ are elements of an one-sided diring 
$(\, R, \, + , \, \xl , \, \xr \,)$ and $\ast \in\{\, \xl ,\, \xr \, \}$.

\begin{definition}\label{def1.0} Let $(\, R, \, + , \, \xl , \, \xr \,)$ be a 
left diring with a left multiplicative bar-unit $e_{\ell}$. The set
$$\hbar ^+ (R)=\{\, x \in R \, | \, e_{\ell}\xl x=0 \,\}$$
is called the {\bf additive halo} of $R$.
\end{definition}

It is clear that the definition of the additive halo $\hbar ^+ (R)$ does not 
dependent on the choice of the left multiplicative bar-unit of $R$. Since 
$$
e_{\ell}\xl x=0\Rightarrow x\xr e_{\ell}=(e_{\ell}\xr x)\xr e_{\ell}=
(e_{\ell}\xl x)\xr e_{\ell}=0\xr e_{\ell}=0,
$$ 
the additive halo $\hbar ^+ (R)$ can be also described as follows
$$\hbar ^+ (R)=\{\, x \in R \, | \, e_{\ell}\xl x=0=x\xr e_{\ell}\,\}.$$

\medskip
The notion of the additive halo is indispensable to rewrite commutative ring 
theory in the context of dirings. The motivation of introducing the notion comes 
from the following facts, which were obtained in our attempt to generalize the 
Lie correspondence between connected linear Lie groups and linear Lie algebras. 

Let $(\, R, \, + ,\, \xl ,\, \xr \,)$ be a diring with a multiplicative 
bar-unit $e$. According to what we did in Section 4.1 of \cite{Liu2}, there are 
three more binary operations
$\xluplus$, $\xruplus$ and $\bullet$ on $R$. Their definitions are as follows:
\begin{eqnarray}
x\xluplus y:&=&x+e\xl y,\nonumber\\
x\xruplus y:&=&x\xr e+ y,\nonumber\\
\label{eq3} x\bullet y: &=& x\xr y +x\xl y - x\xr e\xl y,
\end{eqnarray}
where $x$, $y\in  R$.

One can check that $(\, R, \, \xluplus ,\, \xruplus \, )$ is a 
digroup\footnote{The notion of a digroup we shall use in this paper was introduced 
in Definition 1.1 of \cite{Liu3}. In other words, the left 
inverse of an element $x$ of a digroup may be not equal to the right inverse of $x$.} 
with respect to the bar-unit $0$, and the halo of the digroup is the additive halo 
$\hbar ^+ (R)$ of $R$. The binary operation defined by (\ref{eq3}) is called the 
{\bf Liu product} induced by $+$, $\xl$, $\xr$ and $e$. Since the Liu product is 
associative, a diring $(\, R, \, +,\, \xl ,\xr \,)$ can be regarded as a ring 
$(\, R, \, \bullet \,)$ with the identity $e$. 

\medskip
Let $A$ and $B$ be two subsets of a one-sided diring 
$(\, R, \, + , \, \xl , \, \xr \,)$. We shall use $A\ast B$ to indicate the 
following subset of $R$
$$ A\ast B:=\{\, a\ast b \, |\, a\in A, \, b\in B \,\},$$
where $\ast \in \{\, +, \, \xl ,\, \xr \, \}$. We also use 
$x\equiv y\,  (mod\, A)$ to indicate that $x-y\in A$ for $x$, $y\in R$.

\begin{proposition}\label{pr1.1} Let $(\, R, \, + , \, \xl , \, \xr \,)$ be a 
left diring with a left multiplicative bar-unit $e_{\ell}$.
\begin{description}
\item[(i)] For all $x$, $y\in R$ and $\ast$, 
$\diamond\in \{\, \xl ,\, \xr \, \}$, we have
\begin{equation}\label{eq1} 
x\ast y\equiv x\diamond y \, \left( mod \,\hbar ^+ (R)\right).
\end{equation}
\item[(ii)] $\hbar ^+ (R) \xr R=0=R\xl \hbar ^+ (R)$.
\item[(iii)] $\hbar ^+ (R) \ast R\subseteq \hbar ^+ (R)$ and 
$R\ast \hbar ^+ (R)\subseteq \hbar ^+ (R)$ for 
$\ast\in \{\, \xl ,\, \xr \, \}$.
\item[(iv)] $e_{\ell}+ \hbar ^+ (R)\subseteq \hbar _{\ell}^\times (R)$. 
\item[(v)] If $e\in \hbar ^\times (R)$, then 
$e + \hbar ^+ (R)=\hbar ^\times (R)$. 
\end{description}
\end{proposition}

\medskip
\noindent
{\bf Proof} (i) Since
$$ e_{\ell}\xl (x\ast y\ - x\diamond y)=e_{\ell}\xl x\xl y- e_{\ell}\xl x\xl y=0
$$
for all $x$, $y\in R$ and $\ast$, $\diamond \in \{\, \xl ,\, \xr \, \}$, (i) 
holds.

\medskip
(ii) This part follows from 
$$\hbar ^+ (R) \xr R=\hbar ^+ (R) \xr \left(e_{\ell} \xr R\right)=
\left(\hbar ^+ (R) \xr e_{\ell}\right) \xr R=0\xr R=0$$
and
$$R\xl \hbar ^+ (R) =R\xl \left(e_{\ell} \xr \hbar ^+ (R)\right)=
R\xl \left(e_{\ell} \xl \hbar ^+ (R)\right)=R\xl 0=0. $$

\medskip
(iii) For $\ast\in \{\, \xl ,\, \xr \, \}$, we have
$$ e_{\ell}\xl \left(\hbar ^+ (R) \ast R \right)=
\left(e_{\ell}\xl \hbar ^+ (R)\right) \xl R
=0\xl R=0$$
and 
$$e_{\ell}\xl \left(R\ast \hbar ^+ (R)\right)=e_{\ell}\xl 
\left(R\xl \hbar ^+ (R)\right)=
e_{\ell}\xl 0=0,$$
which imply that (iii) holds.

\medskip
(iv) By (ii), we have
$$ \left( e_{\ell} + \hbar ^+ (R) \right)\xr x=
e_{\ell}\xr x + \hbar ^+ (R)\xr x=x+0=x.$$
Hence, $e_{\ell}+ \hbar ^+ (R)\subseteq \hbar _{\ell}^\times (R)$. 

\medskip
(v) By (iv), we have $e + \hbar ^+ (R)\subseteq \hbar ^\times (R)$.
Conversely, if $\alpha\in \hbar ^\times (R)$, we have
$$e\xl (\alpha -e)=e\xl \alpha -e\xl e=e-e=0,$$
which implies that $\alpha -e\in \hbar ^+ (R)$. Hence, 
$\alpha =e+ (\alpha -e)\in e+ \hbar ^+ (R)$. Thus, 
$\hbar ^\times (R) \subseteq e+ \hbar ^+ (R)$. This proves (v).
\hfill\raisebox{1mm}{\framebox[2mm]{}}

\bigskip
Let $(\, R, \, + , \, \xl , \, \xr \,)$ be a one-sided diring. 
A subgroup $I$ of the additive group $(\, R, \, + \,)$ is called an {\bf ideal} 
of $R$ if 
$$ R\ast I\subseteq I, \qquad I\ast R\subseteq I $$
for $\ast\in \{\, \xl , \, \xr \,\}$. It is clear that if $\hbar ^+(R)\ne 0$, 
then every one-sided diring $R$ always has three distinct ideals: $0$, 
$\hbar ^+(R)$ and $R$ by Proposition~\ref{pr1.1}(iii). 

\begin{definition}\label{def0.4} A one-sided diring $R$ is said to be 
{\bf $3$-simple} if  $\hbar ^+(R)\ne 0$ and $R$ has no ideals other than  $0$, 
$\hbar ^+(R)$ and $R$.
\end{definition}

A one-sided diring $R$ is said to be {\bf $2$-simple} if $R$ has exactly two 
distinct ideals. It is clear that if $R$ is $2$-simple, then $\hbar ^+(R)=0$. 
Hence, the notion of a $2$-simple diring is the same as the notion of a simple ring.

\medskip
Let $I$ be an ideal of a left diring $(\, R, \, + , \, \xl , \, \xr \,)$, and 
let $e_{\ell}$ be a left multiplicative bar-unit of $R$. We define two binary 
operations $\xl$ and $\xr$ on the quotient group
$$\displaystyle\frac{R}{I}:=\left\{\, x+I \, |\, x\in I \, \right\}$$
by
\begin{eqnarray*}
(x+I)\xl (y+I):&=&x\xl y +I,\nonumber\\
(x+I)\xr  (y+I):&=&x\xr y +I,\nonumber
\end{eqnarray*}
where $x$, $y\in  R$. The two binary operations above make the quotient group 
$\displaystyle\frac{R}{I}$ into a left diring with a left multiplicative bar-unit 
$e_{\ell} +I$, which is called the {\bf quotient left diring} of $R$ with 
respect to the ideal $I$.

\medskip
It is clear that if $I$ is an ideal of a left diring $R$ and 
$I\supseteq \hbar ^+(R)$, then the quotient left diring  
$\displaystyle\frac{R}{I}$ is a rng with a left identity.

\medskip
Let $(\, R, \, + , \, \xl , \, \xr \,)$ be a left diring. A subset $S$ of is 
called a {\bf subdiring} of $R$ if $(S, \, +)$ is a subgroup of the additive group 
$(R, \, +)$, $(S, \, \xl, \, \xr)$ is a dimonoid and 
$S\cap \hbar _{\ell}^{\times}(R)\ne\emptyset$. 

\begin{definition}\label{def0.5} Let $R$ and $\bar{R}$ be left dirings. 
 A map $\phi: R\to \bar{R}$ is 
called a {\bf left diring homomorphism} if 
\begin{eqnarray*}
&&\phi (a+b)=\phi (a)+\phi (b),\\
&&\phi (a\ast b)=\phi (a)\ast\phi (b),\\
&&\phi \left(\hbar ^{\times}_{\ell} (R)\right)\cap \hbar ^{\times}_{\ell} 
(\bar{R})\ne \emptyset,
\end{eqnarray*}
where $a$, $b\in R$ and $\ast \in \{ \xl , \, \xr \}$. A bijective left diring 
homomorphism is called a {\bf left diring isomorphism}.
\end{definition}

Let $\phi : R\to \bar{R}$ be a left diring homomorphism from a left diring $R$ 
to a left diring $\bar{R}$. The {\bf kernel} $Ker\phi$ and the {\bf image} 
$Im\phi$ of $\phi$ are defined by
$$ Ker\phi :=\{\, a \,|\, \mbox{$a\in R$ and $\phi (a)=0$} \, \}$$ 
and
$$ Im\phi :=\{\, \phi (a) \,|\, \mbox{$a\in R$} \, \}.$$ 

\medskip
It is clear that $Ker\phi$ is an ideal of the left diring $R$, $Im\phi$ is a 
subdiring of the left diring $\bar{R}$ and
$$\bar{\phi}: a+Ker\phi \mapsto \phi (a) \qquad\mbox{for $a\in R$}$$
is a left diring isomorphism from the quotient diring 
$\displaystyle\frac{R}{Ker\phi}$ to the subdiring of $\bar{R}$.

\medskip
\section{Modules Over One-sided Dirings }

We begin this section with the definition of a left module over a left diring.

\begin{definition}\label{def3.1} Let  $(\, R, \, + , \, \xl , \, \xr \,)$ be a 
left diring with a left multiplicative bar-unit $e$. A {\bf left $R$-module} 
$(\, M, \, \xlodot , \, \xrodot \,)$ is an Abelian group $M$ together with two maps 
$(a, x)\mapsto a\xlodot x$ and $(a, x)\mapsto a\xrodot x$ from $R\times M$ to 
$M$ satisfying the following conditions:
\begin{eqnarray}
\label{eq3.1} a\ast (x+y)&=& a\ast x +a\ast y, \\
\label{eq3.2} (a+b)\ast x&=& a\ast x +x\ast x, \\
\label{eq3.3} (a\xl b)\xlodot x&=& a\xlodot (b\ast x), \\
\label{eq3.4} (a\xr b)\xlodot x&=& a\xrodot (b\xlodot x), \\
\label{eq3.5} (a\diamond b)\xrodot x&=& a\xrodot (b\xrodot x), \\
\label{eq3.6} e\xrodot x&=& x,
\end{eqnarray}
where $a$, $b\in R$, $x$, $y\in M$, $\ast \in \{\xlodot , \, \xrodot \}$ and 
$\diamond \in \{\xl , \, \xr \}$. 
\end{definition}

Let $End(M)$ be the ring of endomorphisms of an Abelian group $M$. If $a$ is an 
element of a left R-module $(\, M, \, \xlodot , \, \xrodot \,)$ over a left 
diring $R$, then both $\Ll _a$ and 
$\Lr _a$ are endomorphisms of $M$, where $\Ll _a$ and $\Lr _a$ are defined by
\begin{equation}\label{eq3.7}
\Ll _a (x):=a\xlodot x, \quad \Lr _a (x):=a\xrodot x \quad\mbox{for $x\in M$.}
\end{equation}
$\Ll _a$ and $\Lr _a$ are called the {\bf left translations} determined by $a$, 
which have been used to study digroups in Chapter 2 of \cite{Liu2}. 

It is easy to check that the two maps $\Ll : a\mapsto \Ll _a$ and 
$\Lr : a\mapsto \Lr _a$ are two group homomorphisms from the additive group 
$(\, R, \, +\, )$ to the additive group 
$(\, End(M), \, +\,)$ and the following are true:
\begin{eqnarray}\label{eq3.8} \Ll _a\Ll _b=\Ll _a\Lr _b=\Ll _{a\xl b}, 
\end{eqnarray}
\begin{eqnarray}
\label{eq3.9} \Lr _a\Lr _b=\Lr _{a\ast b}, \qquad \Lr _a\Ll _b=\Ll _{a\xr b},
\end{eqnarray}
\begin{eqnarray}
\label{eq3.10} \Lr _a\Lr _e=\Lr _a, \qquad \Lr _e=1,
\end{eqnarray}
where $a$, $b\in R$

Conversely, if there are two group homomorphisms $\Ll : a\mapsto \Ll _a$ and 
$\Lr : a\mapsto \Lr _a$ from  $(\, R, \, +\,)$ to $(\,End(M), \, +\,)$ 
satisfying (\ref{eq3.8}), (\ref{eq3.9}) and (\ref{eq3.10}), then the Abelian 
group $M$ becomes a left $R$-module under the module actions defined by (\ref{eq3.7}).

\medskip
Let $R$ be a left diring. A subset $N$ of a left $R$-module 
$(\, M, \, \xlodot , \, \xrodot \,)$
is called a {\bf submodule} of $M$ if the following conditions are satisfied:
\begin{description}
\item[(i)] $N$ is a subgroup of the Abelian group $M$,
\item[(ii)] For all $a\in R$, $x\in N$ and 
$\ast \in \{\xlodot ,\, \xrodot  \}$, $a\ast x\in N$.
\end{description}

\medskip
Let $M$ be a left $R$-module over a left diring $R$. It is clear that both $0$ 
and $M$ are submodules of $M$. Let $K$ and $J$ be two submodules of $M$. We say 
that a submodule $N$ of $M$ is a {\bf proper submodule between $K$ and $J$} 
if $K\ne N$, $N\ne J$ and $K\subseteq N\subseteq J$.

\medskip
Let $e$ be a multiplicative bar-unit of $R$. The {\bf additive halo} 
$\hbar ^+(M)$ of $M$ is defined by
\begin{equation}\label{eq3.11}
\hbar ^+(M): =\{\, x\in M \,|\, e\xlodot x=0\,\}.
\end{equation}

It is immediate that $\hbar ^+(M)$ is a submodule of $M$. Hence, every left 
$R$-module $M$ always has three submodules: $0$, $\hbar ^+(M)$ and $M$.

\medskip
A left diring $(\, R, \, + , \, \xl , \, \xr \,)$ can be regarded as a left 
$R$-module 
$(\, R, \, \xlodot , \, \xrodot \,)$, where $\xlodot :=\xl$ and 
$\xrodot :=\xr$. This module is denoted by $_RR$ and is called the 
{\bf left regular module over $R$}. A submodule of a left regular module 
$_RR$ is called a {\bf left ideal} of $R$.

\begin{proposition}\label{pr3.1} Let $(\, R, \, + , \, \xl , \, \xr \,)$ be a 
left diring with a left multiplicative bar-unit $e$. If 
$(\, M, \, \xlodot , \, \xrodot \,)$ is a left $R$-module, then
\begin{description}
\item[(i)] $R\xlodot \hbar ^+(M)=0$.
\item[(ii)] $\hbar ^+(R)\xrodot M=0$ and 
$\hbar ^+(R)\xlodot M\subseteq \hbar ^+(M)$.
\item[(iii)] $M=(e\xlodot M)\oplus \hbar ^+(M)$, where $\oplus$ denotes the 
direct sum of groups.
\end{description}
\end{proposition}

\medskip
\noindent
{\bf Proof}  (i) For $a\in R$ and $x\in \hbar ^+(M)$, we have
$$ a\xlodot x=a\xlodot (e\xrodot x)=a\xlodot (e\xlodot x)=a\xlodot 0=0$$
by (\ref{eq3.6}) and (\ref{eq3.3}). This proves (i).

\medskip
(ii) For $a\in \hbar ^+(R)$ and $y\in M$, we have
$$a\xrodot y=a\xrodot (e\xrodot y)=(a\xr e)\xrodot y=0\xrodot y=0$$
by (\ref{eq3.6}) and (\ref{eq3.5}). Hence, $\hbar ^+(R)\xrodot M=0$.

Using (\ref{eq3.3}), we have
$$ e\xlodot (a\xlodot y)=(e\xl a)\xlodot y=0\xlodot y=0,$$
which proves that $\hbar ^+(R)\xlodot M\subseteq \hbar ^+(M)$.

\medskip
(iii) For any $z\in M$, we have
\begin{equation}\label{eq3.12}
z=e\xlodot z +(z-e\xlodot z).
\end{equation}

By (\ref{eq3.1}), (\ref{eq3.6}) and (\ref{eq3.3}), we have
\begin{eqnarray*}
&&e\xlodot (z-e\xlodot z)=e\xlodot z -e\xlodot (e\xlodot z)\\
&=&e\xlodot (e\xrodot z) -e\xlodot (e\xlodot z)=0,
\end{eqnarray*} 
which implies that
\begin{equation}\label{eq3.13}
z-e\xlodot z\in \hbar ^+(M) \quad\mbox{for $z\in M$}.
\end{equation}

By (\ref{eq3.12}) and (\ref{eq3.13}), we get
\begin{equation}\label{eq3.14}
M=(e\xlodot M)+ \hbar ^+(M).
\end{equation}

If $e\xlodot u\in (e\xlodot M)\cap \hbar ^+(M)$ with $u\in M$, then
$$
e\xlodot u=e\xlodot (e\xrodot u)=e\xlodot (e\xlodot u)\in R\xlodot \hbar ^+(M)=0
$$
by (i). Hence, we get
\begin{equation}\label{eq3.15}
(e\xlodot M)\cap \hbar ^+(M)=0.
\end{equation}

It follows from (\ref{eq3.14}) and (\ref{eq3.15}) that (iii) holds.
\hfill\raisebox{1mm}{\framebox[2mm]{}}

\bigskip
By Proposition~\ref{pr3.1}(iii), every left multiplicative bar-unit $e$ of a 
left diring $R$ induces a decomposition of a left $R$-module $M$:
\begin{equation}\label{eq3.16}
M=M_0\oplus M_1,
\end{equation}
where
$$M_0:=e\xlodot M, \qquad M_1:=\hbar ^+(M).$$

By (\ref{eq3.16}), every element $x$ of a left $R$-module $M$ can be expressed 
uniquely as 
$$ x=x_0+x_1, \qquad\mbox{ $x_i\in M_i$ for $i=0$, $1$.}$$
$x_0$ and $x_1$ are called the {\bf even component} of $x$ and the 
{\bf odd component} of $x$ induced by $e$, respectively.

\medskip
A useful property of even components is
\begin{equation}\label{eq3.17}
e\xlodot x_0=x_0 \qquad\mbox{for $x_0\in e\xlodot M$.}
\end{equation}

\medskip
Let $M$ and $\bar{M}$ be two left modules over a left diring $R$. A map 
$\phi : M\to \bar{M}$ is called a {\bf $R$-homomorphism} 
(or {\bf module homomorphism}) if
\begin{eqnarray*}
\phi (x+y)&=& \phi (x) +\phi (y),\\
\phi (a\ast x)&=& a\ast \phi (x),
\end{eqnarray*}
for $x$, $y\in M$, $a\in R$ and $\ast\in\{ \xlodot, \, \xrodot \}$. A bijective 
$R$-homomorphism is called a {\bf $R$-isomorphism}. The {\bf kernel} $Ker \phi$ 
and the {\bf image} $Im\phi$ of a $R$-homomorphism $\phi : M\to \bar{M}$ are 
defined by
$$ Ker \phi : =\{\, x\, |\, \mbox{$x\in M$ and $\phi (x)=0$}\,\}$$
and
$$ Im\phi : =\{\, \phi (x)\, |\, \mbox{$x\in M$}\,\}.$$
   
It is easy to check that $Ker \phi$ is a submodule of $M$, $Im\phi$ is a 
submodule of $\bar{M}$ and
\begin{equation}\label{eq3.23}
\phi (\hbar ^+(M))\subseteq \hbar ^+(Im\phi )=\hbar ^+(M)\cap Im\phi.
\end{equation}

\medskip
Let $N$ be a submodule of a left module $(\, M, \, \xlodot , \, \xrodot \,)$ 
over a left diring $R$. Since $a\ast N\subseteq N$ for $a\in R$ and 
$\ast \in \{ \xlodot ,\, \xrodot \}$, we know that
\begin{equation}\label{eq3.24}
a\ast (x+N): = a\ast x +N \qquad\mbox{for $x\in M$}
\end{equation}
is a well defined map from $R\times \left(\displaystyle\frac{M}{N}\right)$ to 
the quotient group $\displaystyle\frac{M}{N}$. One  can check that (\ref{eq3.24}) 
makes $\displaystyle\frac{M}{N}$ into a left $R$-module, which is called the 
{\bf quotient module} of $M$ with respect to the submodule $N$. The additive halo o
f the quotient module $\displaystyle\frac{M}{N}$  is given by
\begin{equation}\label{eq3.25}
\hbar ^+\left( \displaystyle\frac{M}{N}\right)=
\displaystyle\frac{N+\hbar^+(M)}{N}.
\end{equation}

\medskip
\section{$3$-Irreducible Modules}

We now introduce the notion of a $3$-irreducible left module over a left diring.

\begin{definition}\label{def3.2} Let $R$ be a left diring. A left $R$-module 
$M$ is called a {\bf $3$-irreducible module} if $\hbar ^+(M)$ is the unique proper 
submodule between $0$ and $M$. 
\end{definition}

\medskip
Let $R$ be a left diring with a left multiplicative bar-unit $e$ and $M$ a left 
$R$-module. A submodule $N$ of $M$ is said to be {\bf $3$-maximal} if the 
quotient module $\displaystyle\frac{M}{N}$ is $3$-irreducible. By (\ref{eq3.25}), 
a submodule $N$ of $M$ is $3$-maximal if and only if $N+\hbar^+(M)$ is the unique 
proper submodule between $N$ and $M$. A $3$-maximal submodule of the left regular 
module $_RR$ is called a {\bf $3$-maximal left ideal}.

\medskip
The next proposition gives the characterizations of a $3$-irreducible left 
$R$-module.

\begin{proposition}\label{pr3.5} Let $R$ be a left diring with a left 
multiplicative bar-unit $e$. If $(\, M, \, \xlodot , \, \xrodot \,)$ is a left 
$R$-module with $M\ne \hbar ^+(M)$ and 
$\hbar ^+(M)\ne 0$, then the following are equivalent:
\begin{description}
\item[(i)] $M$ is $3$-irreducible.
\item[(ii)] $M=R \xlodot x_0$ for any nonzero element $x_0$ of $e\xlodot M$, 
and $\hbar ^+(M)=$ $R\xrodot x_1$ for any nonzero element $x_1$ of $\hbar ^+(M)$.
\item[(iii)] $M\simeq \displaystyle\frac{R}{I}$ as left $R$-modules, where $I$ 
is a $3$-maximal left ideal of $R$.
\end{description}
\end{proposition}

\medskip
\noindent
{\bf Proof} This is a direct consequence of Definition~\ref{def3.2}.
\hfill\raisebox{1mm}{\framebox[2mm]{}}

\bigskip
Let $\{\, M_{\lambda}' \,|\, \lambda\in \Lambda \,\}$ be a family of left 
modules over a left diring $R$. The {\bf (external) direct sum} 
$\mathring{\displaystyle\bigoplus _{\lambda \in \Lambda}}M_{\lambda}'$ of the 
left $R$-modules $M_{\lambda}'$ is defined by
$$
\mathring{\displaystyle\bigoplus _{\lambda \in \Lambda}}M_{\lambda}':=
\left\{\, f \, \left| \, \begin{array}{c} 
\mbox{$f: \Lambda \to \displaystyle\bigcup _{\lambda \in \Lambda}M_{\lambda}'$ is a map}\\
\mbox{ such that $f(\lambda)\in M_{\lambda}'$  for $\lambda \in \Lambda$ and }\\
\mbox{$supp\, f:=\{\,\lambda \, |\, \mbox{$\lambda \in \Lambda$ and 
$f(\lambda )\ne 0$} \}$ is a finite set}\end{array}\right.\right\}.
$$

For $f$, 
$g\in \mathring{\displaystyle\bigoplus _{\lambda \in \Lambda}}M_{\lambda}'$, 
$a\in R$ and $\ast\in \{\xlodot ,\, \xrodot \}$, we define $f+g$ and 
$a\ast f$ by
\begin{eqnarray}
\label{eq3.27}(f+g)(\lambda):&=& f(\lambda) +g(\lambda),\\
\label{eq3.28}(a\ast f)(\lambda):&=& a\ast f(\lambda),
\end{eqnarray}
where $\lambda \in \Lambda$. It is easy to check that 
$\left( \mathring{\displaystyle\bigoplus _{\lambda \in \Lambda}}
M_{\lambda}', \, \xlodot , \, \xrodot \right)$ is a left $R$-module.

\begin{definition}\label{def3.3} Let $R$ be a left diring. A left $R$-module $M$ 
is said to be {\bf completely $3$-reducible} if $M$ is a direct sum of 
$3$-irreducible left $R$-modules. 
\end{definition}

Let $M$ and $N$ be left modules over a left diring $R$. The set of all 
$R$-homomorphisms from $M$ to $N$ is denoted by $Hom\, _R (M, \, N)$. It is clear that 
$\left(Hom\, _R (M, \, N), \, +, \, 0 \right)$ is an Abelian group, where 
the addition $+$ is defined by
$$
(f +g) (x): =f(x)+g(x) \quad\mbox{for $x\in M$}
$$
and the $R$-homomorphism $0\in Hom\, _R (M, \, N)$ is defined by
$$
0(x):=0 \quad\mbox{for $x\in M$.}
$$
The additive inverse $-f$ of an element $f\in Hom\, _R (M, \, N)$ is given by 
$$
(-f)(x):=- f(x) \quad\mbox{for $x\in M$.}
$$

If $M=N$, the Abelian group
$$
End\, _R M:=Hom\, _R (M, \, M)
$$
is a ring with respect to the associative product $fg$, where $fg$ is defined by
$$
(fg)(x):=f(g(x)) \quad\mbox{for $x\in M$.}
$$

\medskip
The next proposition shows that Schur's Lemma is still true for $3$-irreducible 
left modules over a left diring.

\begin{proposition}\label{pr3.6} Let $R$ be a left diring.
If $M$ and $N$ are $3$-irreducible left $R$-modules, then any $R$-homomorphism 
from $M$ to $N$ is either $0$ or a $R$-isomorphism. In other words, $End\, _R M$ 
is a division ring.
\end{proposition}

\medskip
\noindent
{\bf Proof} Let $f$ be a nonzero $R$-homomorphism from $M$ to $N$. Then 
$Ker f=\hbar ^+(M)$ or $Ker f=0$, and $Im f=\hbar ^+(N)$ or $Im f =N$. 
Hence, we have four possible cases.

\medskip
\underline{{\it Case 1}}: $Ker f=\hbar ^+(M)$ and $Im f=\hbar ^+(N)$, in which 
case, we have
$$
f(e\xlodot x)=e\xlodot f(x)\in e\xlodot \hbar ^+(N)=0 \quad\mbox{for $x\in M$.}
$$

Hence, $e\xlodot M\subseteq Ker f=\hbar ^+(M)$, which is impossible.

\medskip
\underline{{\it Case 2}}: $Ker f=\hbar ^+(M)$ and $Im f=N$, in which case, 
we have
$$
\displaystyle\frac{M}{\hbar ^+(M)}=
\displaystyle\frac{M}{Ker f}\simeq N \quad\mbox{as left $R$-modules.}
$$
Since $\hbar ^+(N)$ is a proper submodule between $0$ and $N$, there is a proper 
submodule between $\hbar ^+(M)$ and $M$, which is impossible.

\medskip
\underline{{\it Case 3}}: $Ker f=0$ and $Im f=\hbar ^+(N)$, in which case, 
$0\ne e\xlodot M\subseteq Ker f=0$, which is impossible.

\medskip
\underline{{\it Case 4}}: $Ker f=0$ and $Im f=\hbar ^+(N)$, in which case, $f$ 
is a $R$--isomorphism.

\hfill\raisebox{1mm}{\framebox[2mm]{}}

\medskip
\section{$3$-Primitivity and $3$-Semi-Primitivity}

Let $(\, R, \, +, \, \xl , \, \xr \,)$ be a left diring. The {\bf annihilator} 
$ann _R M$ of a left $R$-module $(\, M, \, \xlodot , \, \xrodot \,)$ is defined by
\begin{equation}\label{eq4.1}
ann _R M:=\left\{\, a\in R \, |\, a\ast M=0 \quad
\mbox{for $\ast\in \{\xlodot , \, \xrodot \}$} \,\right\}.
\end{equation}

It is clear that $ann _R M$ is an ideal of $R$.

\medskip
A left $R$-module $(\, M, \, \xlodot , \, \xrodot \,)$ is said to be 
{\bf faithful} if $ann _R M=0$. It is easy to check that 
$(\, M, \, \xlodot , \, \xrodot \,)$ is a faithful left 
$\displaystyle\frac{R}{ann _R M}$-module and the module actions are defined by
\begin{eqnarray}
\label{eq4.2}(b+ann _R M)\xlodot x:& =&b\xlodot x, \\
\label{eq4.3}(b+ann _R M)\xrodot x:& =&b\xrodot x,
\end{eqnarray}
where $b\in R$ and $x\in M$.

\begin{definition}\label{def4.1} A left diring $R$ is said to be 
{\bf $3$-primitive} if there is a faithful $3$-irreducible left $R$-module. 
A left diring $R$ is said to be {\bf $3$-semi-primitive} if for any $a\ne 0$ 
in $R$ there exists a $3$-irreducible left $R$-module $M$ such that 
$a\not\in ann _R M$.
\end{definition}

Let $\{\, R_{\lambda} \, |\, \lambda\in \Lambda \, \}$ be a family of left 
dirings indexed by a set $\Lambda$. The set
\begin{equation}\label{eq4.4}
\displaystyle\prod _{\lambda \in \Lambda}R_{\lambda}:=
\left\{\, f \, \left| \, \begin{array}{c} 
\mbox{$f: \Lambda \to \displaystyle\bigcup _{\lambda \in \Lambda}R_{\lambda}$ is a map}\\
\mbox{ such that $f(\lambda)\in R_{\lambda}$  for 
$\lambda \in \Lambda$}\end{array}\right.\right\}.
\end{equation}
is called the {\bf direct product} of the left dirings $R_{\lambda}$ with 
$\lambda\in \Lambda$. For $f$, 
$g\in \displaystyle\prod _{\lambda \in \Lambda}R_{\lambda}$, we define 
$f+g$, $f\xl g$ and $f\xr g$ by
\begin{eqnarray*}
(f+g) (\lambda ) :&=& f (\lambda )+g (\lambda ),\\
(f\xl g) (\lambda ) :&=& f (\lambda )\xl g (\lambda ),\\
(f\xr g) (\lambda ) :&=& f (\lambda )\xr g (\lambda )
\end{eqnarray*}
for all $\lambda \in \Lambda$. Let $0_{\lambda}$ and $e_{\lambda}$ be the zero 
element of $R_{\lambda}$ and a left multiplicative bar-unit of 
$R_{\lambda}$, respectively. We define $0_{\Lambda}$ and $e_{\Lambda}$ by 
$$ 0_{\Lambda} (\lambda ):=0_{\lambda}, \qquad e_{\Lambda} (\lambda ):=
e_{\lambda}\qquad\mbox{for all $\lambda\in \Lambda$}.$$
Then the direct product 
$\left( \displaystyle\prod _{\lambda \in \Lambda}
R_{\lambda},\, + , \, \xl ,\, \xr \right)$ is a left diring, where $0_{\Lambda}$ 
is the zero element of the direct product, and $e_{\Lambda}$ is a left 
multiplicative bar-unit of the direct product. The additive halo and the left 
multiplicative halo of the direct product are given by
\begin{equation}\label{eq4.5}
\hbar ^+\left( \displaystyle\prod _{\lambda \in \Lambda}R_{\lambda}\right)=
\left\{\, f \, | \, \mbox{$f(\lambda)\in \hbar ^+\left( R_{\lambda}\right)$ for 
all $\lambda \in \Lambda$} \right\}
\end{equation}
and
\begin{equation}\label{eq4.6}
\hbar ^\times _{\ell}\left( \displaystyle\prod _{\lambda \in \Lambda}
R_{\lambda}\right)=
\left\{\, f \, | \, 
\mbox{$f(\lambda)\in \hbar ^\times _{\ell}\left( R_{\lambda}\right)$ for all 
$\lambda \in \Lambda$} \right\}.
\end{equation}

For $\lambda \in \Lambda$, the map 
$\pi _{\lambda}: \displaystyle\prod _{\alpha \in \Lambda}R_{\alpha}\to R_{\lambda}$ 
defined by
\begin{equation}\label{eq4.7}
\pi _{\lambda} (f): =f(\lambda) 
\qquad\mbox{for $f\in \displaystyle\prod _{\alpha \in \Lambda}R_{\alpha}$ }
\end{equation}
is a surjective left diring homomorphism. $\pi _{\lambda} $ is called the 
{\bf projection} from
$\displaystyle\prod _{\alpha \in \Lambda}R_{\alpha}$ onto $R_{\lambda}$.

\begin{definition}\label{def4.2} Let 
$\{\, R_{\lambda} \, |\, \lambda\in \Lambda \, \}$ be a family of left dirings 
indexed by a set $\Lambda$. A left diring $R$ 
is called a {\bf subdirect product} of $R_{\lambda}$ with $\lambda\in \Lambda$ if
 there is an injective left diring homomorphism 
$\phi : R\to \displaystyle\prod _{\lambda \in \Lambda}R_{\lambda}$ such that 
$Im (\pi _{\lambda}\phi)=R_{\lambda}$ for all $\lambda\in \Lambda$.
\end{definition}

We now establish two external characterizations of $3$-semi-primitivity.

\begin{proposition}\label{pr4.1} The following conditions on a left diring $R$ 
are equivalent:
\begin{description}
\item[(i)] $R$ is $3$-semi-primitive.
\item[(ii)] There exists a faithful completely $3$-reducible left $R$-module.
\item[(iii)] $R$ is a subdirect product of $3$-primitive left dirings.
\end{description}
\end{proposition}

\medskip
\noindent
{\bf Proof} (i) $\Rightarrow$ (ii): For each $a\ne 0$ in $R$, we have a 
$3$-irreducible left modules $M_a$ such that $a\not\in ann _RM_a$. Form 
$M=\displaystyle\bigoplus _{a \in R\setminus\{0\}}M_a$, which is the direct 
sum of the left modules $M_a$ with $a \in R\setminus\{0\}$. By (\ref{eq3.28}), 
we have
$$ann _RM=\bigcap _{a \in R\setminus\{0\}}ann _RM_a =0.$$

Hence, the direct sum is a faithful completely $3$-reducible left $R$-module.

\medskip
(ii) $\Rightarrow$ (iii): Let $M$ be a faithful completely $3$-reducible left 
$R$-module. Then 
$M=\displaystyle\bigoplus _{\lambda \in \Lambda}M_{\lambda}$ is the direct sum of 
$3$-irreducible left $R$-modules $M_{\lambda}$. Hence, we have
\begin{equation}\label{eq4.8}
0=ann _RM=\bigcap _{\lambda \in \Lambda}ann _RM_{\lambda}.
\end{equation}

Since $ann _RM_{\lambda}$ is an ideal of $R$ for $\lambda \in \Lambda$, we have a
 left diring homomorphism $\phi$ from $R$ to the direct product 
$\displaystyle\prod _{\lambda \in \Lambda}R_{\lambda}$ of the left dirings $R_{\lambda}$, where 
$R_{\lambda}:=\displaystyle\frac{R}{ann _RM_{\lambda}}$ is the quotient left 
diring of $R$ with respect to the ideal $ann _RM_{\lambda}$, and $\phi$ is defined by
\begin{equation}\label{eq4.9}
\phi (a): \lambda \mapsto a+ann _RM_{\lambda} \qquad\mbox{for $a\in R$ and 
$\lambda \in \Lambda$}.
\end{equation}

It follows from (\ref{eq4.8}) and (\ref{eq4.9}) that
\begin{eqnarray*}
&&\phi (a)=0\\
&\Leftrightarrow& a+ann _RM_{\lambda}=ann _RM_{\lambda}
\quad\mbox{for all $\lambda \in \Lambda$}\\
&\Leftrightarrow& a\in ann _RM_{\lambda}
\quad\mbox{for all $\lambda \in \Lambda$}\\
&\Leftrightarrow& a\in \bigcap _{\lambda \in \Lambda}ann _RM_{\lambda}=0,
\end{eqnarray*}
which proves that $\phi$ is injective.

For any $\lambda \in \Lambda$, we have
$$
\pi _{\lambda}(\phi (a))=\phi (a)(\lambda)=a+ann _RM_{\lambda}  
\qquad\mbox{for $a\in R$,}
$$
which implies that 
$$
Im (\pi _{\lambda}\phi)=\displaystyle\frac{R}{ann _RM_{\lambda}}=
R_{\lambda}\qquad\mbox{for all $\lambda \in \Lambda$.}
$$

This proves that $R$ is a subdirect product of left dirings $R_{\lambda}$.

Since $M_{\lambda}$ is a faithful $3$-irreducible left module over 
$\displaystyle\frac{R}{ann _RM_{\lambda}}=R_{\lambda}$ under the module actions 
(\ref{eq4.2}) and (\ref{eq4.3}), $R_{\lambda}$ is a $3$-primitive left diring. 
Therefore, (iii) holds.

\medskip
(iii) $\Rightarrow$ (i): Let $R$ be a subdirect product of the $3$-primitive left
 dirings 
$R_{\lambda}$ with $\lambda\in\Lambda$. Hence, there is an injective left diring 
homomorphism $\phi : R\to \displaystyle\prod _{\lambda \in \Lambda}R_{\lambda}$. 

Let $M_{\lambda}$ be a faithful $3$-irreducible left $R_{\lambda}$-module, where 
$\lambda\in\Lambda$. For $x _{\lambda}\in M_{\lambda}$, we define
\begin{equation}\label{eq4.10}
a\ast x _{\lambda}:=(\pi _{\lambda}\phi)(a) \ast x _{\lambda}=
\phi (a)(\lambda)\ast x _{\lambda},
\end{equation}
where $a\in R$, $\ast\in \{ \xlodot, \, \xrodot \}$, and 
$\pi _{\lambda} : \displaystyle\prod _{\alpha \in \Lambda}R_{\alpha}\to 
R_{\lambda}$ is the projection defined by (\ref{eq4.7}). It is 
clear that $M_{\lambda}$ becomes a left $R$-module under (\ref{eq4.10}).

Note that 
$\phi \left(\hbar ^\times _{\ell}(R)\right)\bigcap \hbar ^
\times _{\ell}\left( \displaystyle\prod _{\lambda \in \Lambda}R_{\lambda}\right)
\ne \emptyset$. Hence, we have $e\in \hbar ^\times _{\ell}(R)$ such that 
$\phi (e)\in \hbar ^\times _{\ell}
\left( \displaystyle\prod _{\lambda \in \Lambda}R_{\lambda}\right)$. 
By (\ref{eq4.6}), $\phi (e) (\lambda)\in \hbar ^\times _{\ell}(R _{\lambda})$. 
For $x_{\lambda}\in M_{\lambda}$, we have
$$
e\xlodot x_{\lambda}=0 \Leftrightarrow \phi (e) (\lambda) \xlodot x_{\lambda}=0
$$
by (\ref{eq4.10}). This proves that the additive halo of the $R$-module 
$M_{\lambda}$ is equal to the additive halo of the $R_{\lambda}$-module $M_{\lambda}$.

Since $Im (\pi _{\lambda}\phi )=R_{\lambda}$ for $\lambda \in \Lambda$, a 
subgroup $N$ of 
$( M_{\lambda}, \, +)$ is a left $R_{\lambda}$-submodule of $M_{\lambda}$ if and 
only if $N$ is a 
left $R$-submodule of $M_{\lambda}$ by (\ref{eq4.10}). This proves that 
$M_{\lambda}$ is a $3$-irreducible left $R$-module under (\ref{eq4.10}).

We now consider the direct sum 
$M=\displaystyle\bigoplus _{\lambda \in \Lambda}M_{\lambda}$ of the 
$3$-irreducible left $R$-modules $M_{\lambda}$. Using 
(\ref{eq4.10}) and the fact that $M_{\lambda}$  is a faithful $R_{\lambda}$-module, 
we get
\begin{eqnarray*}
&&a\in ann _RM \\
&\Leftrightarrow& \pi _{\lambda}(\phi (a))\in ann _{R_{\lambda}}M_{\lambda}=
0\quad\mbox{for $\lambda\in \Lambda$}\\
&\Leftrightarrow& \phi (a)\in \bigcap _{\lambda \in \Lambda} Ker \pi _{\lambda} =
0\\
&\Leftrightarrow& a=0.
\end{eqnarray*}

Hence, $M=\displaystyle\bigoplus _{\lambda \in \Lambda}M_{\lambda}$ is a 
faithful completely $3$-irreducible $R$-module. In other words, $R$ is 
$3$-semi-primitive.
\hfill\raisebox{1mm}{\framebox[2mm]{}}

\medskip
\section{$3$-Primitive Ideals }

Let $I$ be a left ideal of a left diring $(\, R, \, +, \, \xl , \, \xr \,)$, 
We define
\begin{equation}\label{eq4.11}
(I: R):= \{ \, a\in R \, |\, 
\mbox{$a\xl R\subseteq I$ and $a\xr R\subseteq I$} \, \}.
\end{equation}

After regarding $I$ as a submodule of the left regular $R$-module $_RR$, the 
annihilator of the quotient $R$-module $\displaystyle\frac{R}{I}$ is $(I: R)$. Thus, 
we know that
\begin{equation}\label{eq4.12}
(I: R)= ann _R\left( \displaystyle\frac{R}{I}\right)
\end{equation}
is an ideal of $R$. If $K$ is an ideal of $R$ contained in $I$, then 
$K\diamond R\subseteq K\subseteq I$ for $\diamond \in \{\xl , \, \xr \}$. Hence, 
$K\subseteq (I: R)$ by (\ref{eq4.11}). This proves that 
\begin{equation}\label{eq4.13}
\mbox{$(I: R)$ is the largest ideal of $R$ contained in $I$,}
\end{equation}
where $I$ is a left ideal of a left diring $ R$.

Note that
$$
\mbox{$R$ is a diring $\Rightarrow$ $(I: R)\subseteq I$.}
$$

\medskip
Let $(\, M, \, \xlodot , \, \xrodot \,)$ be a left module over a 
left diring $R$. If $H$ is an ideal of $R$ and $H\subseteq ann _RM$, then $M$ is 
a left module over the quotient left diring $\bar{R}:=\displaystyle\frac{R}{H}$ 
under the following module actions:
\begin{equation}\label{eq4.14}
(a+H)\ast x:= a\ast x,
\end{equation}
where $a\in R$, $x\in M$ and $\ast \in \{\xlodot ,\, \xrodot \}$. It is clear 
that
\begin{eqnarray}
&&\mbox{$M$ is a $3$-irreducible left $R$-module}\nonumber\\
\label{eq4.15}&\Leftrightarrow& \mbox{ $M$ is a $3$-irreducible left 
$\bar{R}$-module}
\end{eqnarray}
and
\begin{equation}\label{eq4.16}
ann _{\bar{R}}M=\displaystyle\frac{ann _RM}{H}.
\end{equation}

\begin{definition}\label{def4.3} 
An ideal $H$ of a left diring $R$ is called a {\bf $3$-primitive ideal} if the 
quotient left diring $\displaystyle\frac{R}{H}$ is  a $3$-primitive left diring. 
\end{definition}

\begin{proposition}\label{pr4.4} Let $H$ be an ideal of a left diring $R$. Then 
$H$ is $3$-primitive if and only if $H=(I:R)$ for some $3$-maximal left ideal $I$ 
of $R$.
\end{proposition}

\medskip
\noindent
{\bf Proof} If $H$ is a $3$-primitive ideal, then there exists a $3$-irreducible 
$\bar{R}:=\displaystyle\frac{R}{H}$-module $_{\bar{R}}M$ such that 
$ann _{\bar{R}}M=\{H\}$.

It is clear that $M$ becomes a left $R$-module $_RM$ under the following module 
actions:
\begin{equation}\label{eq4.17}
a\ast x:=(a+H)\ast x,
\end{equation}
where $a\in R$, $x\in M$ and $\ast\in\{ \xlodot, \, \xrodot \}$. Since
$$
a\in ann _RM\Leftrightarrow a+H\in ann _{\bar{R}}M=\{H\}\Leftrightarrow a+H=
H\Leftrightarrow a\in H,
$$
we have $ann _RM=H$. By (\ref{eq4.15}), $M$ is also $3$-irreducible as a left 
$R$-module. Using Proposition~\ref{pr3.5}(iii), $M\simeq \displaystyle\frac{R}{I}$ 
as left $R$-modules, where $I$ is a $3$-maximal left ideal of $R$. Thus, we get
$$
H=ann _RM=ann _R\left(\displaystyle\frac{R}{I}\right)=(I:R).
$$

Conversely, if $H=(I:R)$ for a $3$-maximal left ideal of $R$, then 
$M\simeq \displaystyle\frac{R}{I}$ is a $3$-irreducible left $R$-module such that
$$ann _RM=ann _R\left(\displaystyle\frac{R}{I}\right)=(I:R)=H.$$
Using (\ref{eq4.14}), (\ref{eq4.15}) and (\ref{eq4.16}), $M$ is a faithful 
$3$-irreducible left module over the quotient left diring $\displaystyle\frac{R}{H}$. 
This prove that $\displaystyle\frac{R}{H}$ is a $3$-primitive left diring. Hence, $H$ 
is a $3$-primitive ideal
\hfill\raisebox{1mm}{\framebox[2mm]{}}

\bigskip
Let $R$ be a left diring. The intersection of the annihilators of all 
$3$-irreducible left $R$-modules is called the {\bf $3$-radical} of $R$ and is 
denoted by $rad_3R$. Since
\begin{equation}\label{eq4.18}
rad _3R=\bigcap _{\begin{array}{c}\mbox{$M$ runs over all}\\
\mbox{$3$-irreducible left $R$-module}\end{array}} ann _RM
\end{equation}
and $ann _RM$ is an ideal of $R$, $rad _3R$ is an ideal of $R$.

\begin{proposition}\label{pr4.5} If $R$ is a left diring, then $rad _3 R$ is the 
intersection of the $3$-primitive ideals of $R$.
\end{proposition}

\medskip
\noindent
{\bf Proof} By Proposition~\ref{pr3.5}(iii), (\ref{eq4.12}) and (\ref{eq4.18}), 
we have
$$
rad _3R=\bigcap _{\begin{array}{c}\mbox{$I$ runs over all}\\\mbox{$3$-maximal 
ideal of $R$}\end{array}}(I:R),
$$
which can be written as
\begin{equation}\label{eq4.19}
rad _3R=\bigcap _{\begin{array}{c}\mbox{$H$ runs over all}\\
\mbox{$3$-primitive ideal of $R$}\end{array}}H
\end{equation}
by Proposition~\ref{pr4.4}.
\hfill\raisebox{1mm}{\framebox[2mm]{}}

\begin{proposition}\label{pr4.6} Let $R$ be a nonzero left diring.
\begin{description}
\item[(i)] $R$ is $3$-semi-primitive if and only if $rad _3R=0$.
\item[(ii)] $rad_3\left(\displaystyle\frac{R}{rad _3R}\right)=0$.
\end{description}
\end{proposition}

\medskip
\noindent
{\bf Proof} (i) If $rad _3R=0$, then 
$$
\bigcap _{\begin{array}{c}\mbox{$H$ runs over all }\\
\mbox{$3$-primitive ideal of $R$}\end{array}}H=0
$$
by (\ref{eq4.19}). Hence, $R$ is a subdirect product of the $3$-primitive left 
dirings $\displaystyle\frac{R}{H}$, where $H$ runs over the $3$-primitive ideals 
of $R$. It follows from Proposition~\ref{pr4.1} that $R$ is $3$-semi-primitive.

Conversely, if $R$ is $3$-semi-primitive, then $R$ is a subdirect product of the 
$3$-primitive left dirings $R_{\lambda}$ with $\lambda\in\Lambda$. Hence, there 
is an injective left diring homomorphism 
$\phi : R\to \displaystyle\prod _{\lambda \in \Lambda}R_{\lambda}$ such that 
$Im (\pi _{\lambda}\phi) =R_{\lambda}$. Thus, we have 
$\displaystyle\frac{R}{Ker (\pi _{\lambda}\phi)}\simeq R_{\lambda}$ as left 
dirings. This proves that $Ker (\pi _{\lambda}\phi)$ is a $3$-primitive ideal of $R$. 
If $a\in \displaystyle\bigcap _{\lambda\in\Lambda} Ker (\pi _{\lambda}\phi)$, 
then
$$
0=\pi _{\lambda}(\phi (a))=\phi (a) (\lambda) 
\qquad\mbox{for $\lambda\in \Lambda$},
$$
which proves that $\phi (a)$ is the zero element of the diring 
$\displaystyle\prod _{\lambda \in \Lambda}R_{\lambda}$. Since $\phi$ is injective, 
$a=0$. It follows from (\ref{eq4.19}) that
$$
0=\displaystyle\bigcap _{\lambda\in\Lambda} Ker (\pi _{\lambda}\phi)\supseteq 
\bigcap _{\begin{array}{c}\mbox{$H$ runs over all}\\
\mbox{$3$-primitive ideal of $R$}\end{array}}H
=rad _3R.
$$

Hence, we get $rad _3R=0$.

\medskip
(ii) $\bar{H}$ is an ideal of $\bar{R}:=\displaystyle\frac{R}{rad _3R}$ if and 
only if 
$\bar{H}=\displaystyle\frac{H}{rad _3R}$ for some ideal $H$ of $R$ containing 
$rad _3R$. Moreover, we have
$$
\displaystyle\frac{R}{H}\simeq\displaystyle\frac{\displaystyle\frac{R}{rad _3R}}
{\displaystyle\frac{H}{rad _3R}}= 
\displaystyle\frac{\bar{R}}{\bar{H}}\quad\mbox{as left dirings.}
$$

By Proposition~\ref{pr4.5}, every $3$-primitive ideal $H$ of $R$ contains 
$rad _3R$. Hence, we have
\begin{eqnarray}
rad _3\bar{R}&=&\bigcap _{\begin{array}{c}\mbox{$\bar{H}$ runs over all}\\ 
\mbox{$3$-primitive ideals of $\bar{R}$}\end{array}}\bar{H}\nonumber\\
\label{eq4.20}&=&
\bigcap _{\begin{array}{c}\mbox{$H$ runs over all}\\ 
\mbox{$3$-primitive ideals of $R$}\end{array}}\left(\displaystyle\frac{H}{rad _3R}\right).
\end{eqnarray}

If $a+rad _3R\in rad _3\bar{R}$, then $a\in H$ for any $3$-primitive ideal $H$ of
 $R$ by (\ref{eq4.20}). Hence, we get
$$
a\in \bigcap _{\begin{array}{c}\mbox{$H$ runs over all}\\ 
\mbox{$3$-primitive ideals of $R$}\end{array}} H=rad _3R.
$$
Thus, $a+rad _3R=rad _3R$ is the zero element of $\bar{R}$. This proves (ii).
\hfill\raisebox{1mm}{\framebox[2mm]{}}

\end{document}